\documentclass[11pt,a4paper]{article}
\usepackage[english,french]{babel}
\usepackage[latin1]{inputenc}
\usepackage[centertags]{amsmath}
\usepackage{amsfonts}
\usepackage{amssymb}
\usepackage{amsthm}
\usepackage{newlfont}

\usepackage{color}      
\usepackage{epsfig}

\input{psfig}

\usepackage{graphicx}

\usepackage{float}

\newtheorem{theo}{Theorem}
\newtheorem{prop}{Proposition}

\newtheorem{lem}{Lemma}

\newcommand{\complex}{\mathbb C}
\newcommand{\integer}{\mathbb N}
\newcommand{\rinteger}{\mathbb Z}
\newcommand{\real}{\mathbb R}

\newcommand{\pochamer}[3]{\left(#1;#2\right)_{#3}}

\newcommand{\eq}{\mathbb{E}_q}

\newcommand{\qhyper}[4]{\vphantom{}_2 \phi_1 \left(#1, #2; #3 ; #4 \right)}
\newcommand{\qhyperc}[3]{\zeta(#1,#2;#3)}
\newcommand{\qhyperabc}[3]{\xi(#1,#2;#3)}
\newcommand{\qhypermatrice}[4]{A(#1,#2;#3;#4)}

\newcommand{\pmatrice}[4]{\begin{pmatrix}#1&#2\\#3&#4 \end{pmatrix}}
\newcommand{\vect}[2]{\left(\begin{array}{*{100}c} #1 \\ #2 \end{array}\right)}

\newcommand{\qcar}[1]{e_{q,#1}}

\newcommand{\fz}[4]{F^{(0)}(#1,#2;#3;#4)}
\newcommand{\finf}[4]{F^{(\infty)}(#1,#2;#3;#4)}

\newcommand{\yz}[4]{Y^{(0)}(#1,#2;#3;#4)}
\newcommand{\yinf}[4]{Y^{(\infty)}(#1,#2;#3;#4)}

\newcommand{\jz}[1]{J^{(0)}(#1)}
\newcommand{\jinf}[2]{J^{(\infty)}(#1,#2)}

\newcommand{\un}{u}
\newcommand{\deux}{v}
\newcommand{\trois}{w}
\newcommand{\quatre}{y}

\newcommand{\Sl}{\text{Sl}}
\newcommand{\Gl}{\text{Gl}}

\newcounter{nbrecasun}
\addtocounter{nbrecasun}{5}

\newcommand{\condun}[1]{\noindent (\textbf{Case #1})}

\newcommand{\conddeux}[1]{\addtocounter{nbrecasun}{#1} \noindent (\textbf{Case
\arabic{nbrecasun}})
\addtocounter{nbrecasun}{-#1}}

\newcounter{nbrecasdeux}
\addtocounter{nbrecasdeux}{8}

\newcommand{\condtrois}[1]{\addtocounter{nbrecasdeux}{#1} \noindent (\textbf{Case
\arabic{nbrecasdeux}})
\addtocounter{nbrecasdeux}{-#1}}

\newcommand{\bon}{\Omega}

\newcommand{\und}[1]{\underline{#1}}
\begin{document}
\selectlanguage{english}
\begin{center}
\begin{huge}Galois groups of the basic hypergeometric equations \footnote{A definitive version of this paper will appear in \textit{Pacific Journal of Mathematics}.} \end{huge}\\
$_{}$\\
\begin{LARGE}by Julien Roques
 \end{LARGE}\\
$_{}$\\
 20th of August 2007\\
\end{center}
$_{}$\\
\noindent \textbf{Abstract}. \textit{In this paper we compute the
Galois groups of basic
hypergeometric equations.}\\

In this paper $q$ is a complex number such that $0<|q|<1$.

\section{Basic hypergeometric series and equations}
The theory of hypergeometric functions and equations dates back at
least as far as Gauss. It has long been and is still an integral
part of the mathematical literature. In particular, the Galois
theory of (generalized) hypergeometric equations attracted the
attention of many authors. For this issue, we refer the reader to
\cite{beukersheckman,beukersbrownheck,katzexpsum} and to the
references therein. We also single out the papers
\cite{duvalmitschi,mitschihyperconf}, devoted to the calculation
of some Galois groups by means of a density theorem (Ramis
theorem).

In this paper we focus our attention on the Galois theory of the
basic hypergeometric equations, the later being natural
$q$-analogues of the hypergeometric equations.

The \textit{basic hypergeometric series}
$\phi(z)=\qhyper{a}{b}{c}{z}$ with parameters $(a,b,c) \in
(\complex^*) ^3$ defined by :
\begin{eqnarray*}\qhyper{a}{b}{c}{z}&=&\sum_{n=0}^{+\infty}
\frac{\pochamer{a,b}{q}{n} }{\pochamer{c,q}{q}{n} } z^n\\
&=& \sum_{n=0}^{+\infty} \frac{(1-a)(1-aq)\cdots
(1-aq^{n-1})(1-b)(1-bq)\cdots (1-bq^{n-1})}{(1-q)(1-q^2)\cdots
(1-q^{n})(1-c)(1-cq)\cdots (1-cq^{n-1})} z^n
\end{eqnarray*}
was first introduced by Heine and was later generalized by
Ramanujan. As regards functional equations, the basic
hypergeometric series provides us with a solution of the following
second order $q$-difference equation, called the \textit{basic
hypergeometric equation} with parameters $(a,b,c)$ :
\begin{equation}\label{equa hypergeo} \phi(q^2z) - \frac{(a+b)z-(1+c/q)}{abz-c/q} \phi(qz) +  \frac{z-1}{abz-c/q}\phi(z)=0.\end{equation}
This functional equation is equivalent to a functional system.
Indeed, with the notations :
$$\lambda(a,b;c;z)= \frac{(a+b)z-(1+c/q)}{abz-c/q}, \ \ \ \
\mu(a,b;c;z)=\frac{z-1}{abz-c/q},$$ a function $\phi$ is solution
of $(\ref{equa hypergeo})$ if and only if the vector
$\Phi(z)=\binom{\phi(z)}{\phi(qz)}$ satisfies the functional
system : \begin{equation}\label{syst hypergeo} \Phi(qz) =
\qhypermatrice{a}{b}{c}{z} \Phi(z) \end{equation} with :
$$\qhypermatrice{a}{b}{c}{z}=\left( \begin{array}{cc} 0&1\\ -\mu(a,b;c;z) & \lambda(a,b;c;z) \end{array}\right).$$
The present paper focuses on the calculation of the Galois group
of the $q$-difference equation (\ref{equa hypergeo}) or,
equivalently, that of the $q$-difference system (\ref{syst
hypergeo}). A number of authors have developed $q$-difference
Galois theories over the past years, among whom Franke
\cite{frankepvdifference}, Etingof \cite{etingofgalois}, Van der
Put and Singer \cite{psgaloistheory}, Van der Put and Reversat
\cite{prev}, Chatzidakis and Hrushovski \cite{chatzihru}, Sauloy
\cite{sauloyqgaloisfuchs}, Andr\'e \cite{andrenoncomm}, etc. The
exact relations between the existing Galois theories for
$q$-difference equations are partially understood. For this
question, we refer the reader to \cite{chatziharsing}, and also to
our Remark section \ref{constr}.

In this paper we follow the approach of Sauloy (initiated by
Etingof in the regular case). Our method for computing the Galois
groups of the basic hypergeometric equations is based on a
$q$-analogue of Schlesinger's density theorem stated and
established in \cite{sauloyqgaloisfuchs}. Note that some of these
groups were previously computed by Hendriks in
\cite{hendriksqhyper} using a radically different method
(actually, the author dealt with the Galois groups defined by Van
der Put and Singer, but these do coincide with those defined by
Sauloy : see our Remark section \ref{constr}). On a related topic, 
we also point out the appendix of \cite{divizio} which contains the 
$q$-analogue of Schwarz's list. 

The paper is organized as follows. In a first part, we give a
brief overview of some results from \cite{sauloyqgaloisfuchs}. In
a second part, we compute the Galois groups of the basic
hypergeometric equations in all non-resonant (but possibly
logarithmic) cases.

\section{Galois theory for regular singular $q$-difference equations}

Using analytic tools together with Tannakian duality, Sauloy
developed in \cite{sauloyqgaloisfuchs} a Galois theory for regular
singular $q$-difference systems. In this section, we shall first
recall the principal notions used in \cite{sauloyqgaloisfuchs},
mainly the Birkhoff matrix and the twisted Birkhoff matrix. Then
we shall explain briefly that this lead to a Galois theory for
regular singular $q$-difference systems. Last, we shall state a
density theorem for these Galois groups, which will be of main
importance in our calculations.

\subsection{Basic notions}\label{section the basic objects}

Let us consider $A \in \Gl_n(\complex(\{z\}))$. Following Sauloy
in \cite{sauloyqgaloisfuchs}, the $q$-difference system :
\begin{equation}\label{general q diff system}Y(qz)=A(z)Y(z)\end{equation}
is said to be \textit{Fuchsian} at $0$ if $A$ is holomorphic at
$0$ and if $A(0)\in \Gl_n(\complex)$. Such a system is
non-resonant at $0$ if, in addition, $Sp(A(0)) \cap
q^{\rinteger^*} Sp(A(0))=\emptyset$. Last we say that the above
$q$-difference system is \textit{regular singular} at $0$ if there
exists $R^{(0)}\in \Gl_n(\complex(\{z\}))$ such that the
$q$-difference system defined by
$(R^{(0)}(qz))^{-1}A(z)R^{(0)}(z)$ is Fuchsian at $0$. We have
similar notions at $\infty$ using the change of variable $z
\leftarrow 1/z$.

In the case of a global system, that is $A\in \Gl_n(\complex(z))$,
we will use the following terminology. If $A\in
\Gl_n(\complex(z))$, then the system (\ref{general q diff system})
is called \textit{Fuchsian} (resp. \textit{Fuchsian and non-resonnant}, \textit{regular
singular}) if it is Fuchsian (resp. Fuchsian and non-resonnant,
regular singular) both at $0$ and at $\infty$.

For instance, the basic hypergeometric system (\ref{syst hypergeo}) is Fuchsian.\\

\textit{Local fundamental system of solutions at $0$}. Suppose
that (\ref{general q diff system}) is Fuchsian and non-resonant at
$0$ and consider $J^{(0)}$ a Jordan normal form of $A(0)$.
According to \cite{sauloyqgaloisfuchs} there exists $F^{(0)} \in
\Gl_n(\complex\{z\})$ such that :
\begin{equation} \label{transfo jauge}
F^{(0)}(qz)J^{(0)}=A(z)F^{(0)}(z).
\end{equation}
Therefore, if $e^{(0)}_{J^{(0)}}$ denotes a fundamental system of
solutions of the $q$-difference system with constant coefficients
$X(qz)=J^{(0)}X(z)$, the matrix-valued function
$Y^{(0)}=F^{(0)}e^{(0)}_{J^{(0)}}$ is a fundamental system of
solutions of (\ref{general q diff system}). We are going to
describe a possible choice for $e^{(0)}_{J^{(0)}}$. We denote by
$\theta_q$ the Jacobi theta function defined by
$\theta_q(z)=\pochamer{q}{q}{\infty}\pochamer{z}{q}{\infty}\pochamer{q/z}{q}{\infty}$.
This is a meromorphic function over $\complex^*$ whose zeros are
simple and located on the discrete logarithmic spiral
$q^\rinteger$. Moreover, the functional equation
$\theta_q(qz)=-z^{-1}\theta_q(z)$ holds. Now we introduce, for all
$\lambda \in \complex^*$ such that $|q| \leq |\lambda| < 1$, the
$q$-character
$e^{(0)}_{\lambda}=\frac{\theta_q}{\theta_{q,\lambda}}$ with
$\theta_{q,\lambda}(z)=\theta_q(\lambda z)$ and we extend this
definition to an arbitrary non-zero complex number $\lambda \in
\complex^*$ requiring the equality
$e^{(0)}_{q\lambda}=ze^{(0)}_{\lambda}$. If
$D=P\text{diag}(\lambda_1,...,\lambda_n)P^{-1}$ is a semisimple
matrix then we set
$e^{(0)}_{D}:=P\text{diag}(e^{(0)}_{\lambda_1},...,e^{(0)}_{\lambda_n})P^{-1}$.
It is easily seen that this does not depend on the chosen
diagonalization. Furthermore, consider
$\ell_q(z)=-z\frac{\theta_q'(z)}{\theta_q(z)}$ and, if $U$ is a
unipotent matrix, $e^{(0)}_{U}=\Sigma_{k=0}^n \ell_q^{(k)}
(U-I_n)^k$ with $\ell_q^{(k)}=\binom{\ell_q}{k}$. If
$J^{(0)}=D^{(0)}U^{(0)}$ is the multiplicative Dunford
decomposition of $J^{(0)}$, with $D^{(0)}$ semi-simple and
$U^{(0)}$ unipotent, we set
$e^{(0)}_{J^{(0)}}=e^{(0)}_{D^{(0)}}e^{(0)}_{U^{(0)}}$.\\

\textit{Local fundamental system of solutions at $\infty$}. Using
the variable change $z \leftarrow 1/z$, we have a similar
construction at $\infty$. The corresponding fundamental system of
solutions is denoted by
$Y^{(\infty)}=F^{(\infty)}e^{(\infty)}_{J^{(\infty)}}$.\\

Throughout this section we assume that the system (\ref{general q
diff system}) is global and that it is Fuchsian
and non-resonant.\\

\textit{Birkhoff matrix}. The linear relations between the two
fundamental systems of solutions introduced above are given by the
Birkhoff matrix (also called connection matrix)
$P=(Y^{(\infty)})^{-1}Y^{(0)}$. Its entries are elliptic functions
\textit{i.e.} meromorphic functions over the elliptic curve
$\eq=\complex^* / q^\rinteger$. \\

\textit{Twisted Birkhoff matrix}. In order to describe a
Zariki-dense set of generators of the Galois group associated to
the system (\ref{general q diff system}), we introduce a
``twisted" connection matrix. According to
\cite{sauloyqgaloisfuchs}, we choose for all $z \in \complex^*$ a
group endomorphism $g_z$ of $\complex^*$ sending $q$ to $z$.
Before giving an explicit example, we have to introduce more
notations. Let us, for any fixed $\tau \in \complex$ such that
$q=e^{-2\pi i \tau}$, write $q^y=e^{-2\pi i \tau y}$ for all $y
\in \complex$. We also define the (non continuous) function
$\log_q$ on the whole punctured complex plane $\complex^*$ by
$\log_q(q^y)=y$ if $y\in \complex^* \setminus \real^+$ and we
require that its discontinuity is located just before the cut
(that is $\real^+$) when turning counterclockwise around $0$. We
can now give an explicit example of endomorphism $g_z$ namely the
function $g_z: \complex^*=\mathbb{U} \times q^\real \rightarrow
\complex^*$ sending $uq^\omega$ to
$g_z(uq^\omega)=z^{\omega}=e^{-2\pi i \tau \log_q(z) \omega}$ for
$(u,\omega)\in \mathbb{U} \times \real$, where $\mathbb{U}\subset
\complex$ is the unit circle.

Then we set, for all $z$ in $\complex^*$,
$\psi_z^{(0)}(\lambda)=\frac{\qcar{\lambda}(z)}{g_z(\lambda)}$ and
we define  $\psi_z^{(0)}\left(D^{(0)}\right)$, the \textit{twisted
factor} at $0$, by
$\psi_z^{(0)}\left(D^{(0)}\right)=P\text{diag}(\psi_z^{(0)}(\lambda_1),...,\psi_z^{(0)}(\lambda_n))P^{-1}$
with $D^{(0)}=P\text{diag}(\lambda_1,...,\lambda_n)P^{-1}$. We
have a similar construction at $\infty$ by using the variable
change $z \leftarrow 1/z$. The corresponding twisting factor is
denoted by $\psi_z^{(\infty)}(J^{(\infty)})$.

Finally, the twisted connection matrix $\breve{P}(z)$ is :
\begin{eqnarray*}
\breve{P}(z)&=&\psi_z^{(\infty)}\left(D^{(\infty)}\right)P(z)\psi_z^{(0)}\left(D^{(0)}\right)^{-1}.
\end{eqnarray*}

\subsection{Definition of the Galois groups}\label{constr}

The definition of the Galois groups of regular singular
$q$-difference systems given by Sauloy in
\cite{sauloyqgaloisfuchs} is somewhat technical and long. Here we
do no more than describe the underlying idea.\\

\textit{(Global) Galois group.} Let us denote by $\mathcal{E}$ the
category of regular singular $q$-difference systems with
coefficients in $\complex(z)$ (so, the base field is
$\complex(z)$; the difference field is $(\complex(z),f(z) \mapsto
f(qz))$). This category is naturally equipped with a tensor
product $\otimes$ such that $(\mathcal{E},\otimes)$ satisfies all
the axioms defining a Tannakian category over $\complex$ except
the existence of a \textit{fiber functor} which is not obvious.
This problem can be overcome using an analogue of the
Riemann-Hilbert correspondance.

The Riemann-Hilbert correspondance for regular singular
$q$-difference systems entails that $\mathcal{E}$ is equivalent to
the category $\mathcal{C}$ of connection triples whose objects are
triples $(A^{(0)},P,A^{(\infty)}) \in \Gl_n(\complex) \times
\Gl_n(\mathcal{M}(\eq)) \times \Gl_n(\complex)$ (we refer to
\cite{sauloyqgaloisfuchs} for the complete definition of
$\mathcal{C}$). Furthermore $\mathcal{C}$ can be endowed with a
tensor product $\und \otimes$ making the above equivalence of
categories compatible with the tensor products. Let us emphasize
that $\und \otimes$ is not the usual tensor product for matrices.
Indeed some twisting factors appear because of the bad
multiplicative properties of the $q$-characters $e_{q,c}$ : in
general $e_{q,c}e_{q,d}\neq e_{q,cd}$.

The category $\mathcal{C}$ allows us to define a Galois group :
$\mathcal{C}$ is a Tannakian category over $\complex$. The functor
$\omega_0$ from $\mathcal{C}$ to $Vect_\complex$ sending an object
$(A^{(0)},P,A^{(\infty)})$ to the underlying vector space
$\complex^n$ on which $A ^{(0)}$ acts is a fiber functor. Let us
remark that there is a similar fiber functor $\omega_\infty$ at
$\infty$. Following the general formalism of the theory of
Tannakian categories (see \cite{deligne}), the \textit{absolute
Galois group} of $\mathcal{C}$ (or, using the above equivalence of
categories, of $\mathcal E$) is defined as the pro-algebraic group
$Aut^{\und \otimes}(\omega_0)$ and the \textit{global Galois group
of an object $\chi$ }of $\mathcal{C}$ (or, using the above
equivalence of categories, of an object of $\mathcal{E}$) is the
complex linear algebraic group $Aut^{\und
\otimes}(\omega_{0|\langle \chi \rangle})$ where $\langle \chi
\rangle$ denotes the Tannakian subcategory of $\mathcal{C}$
generated by $\chi$. For the sake of simplicity, we will often
call $Aut^{\und \otimes}(\omega_{0|\langle \chi \rangle})$ the
\textit{Galois group} of $\chi$
(or, using the above equivalence of categories, of the corresponding object of $\mathcal{E}$). \\

\textit{Local Galois groups.} Let us point out that notions of
local Galois groups at $0$ and at $\infty$ are also available
(here the difference fields are respectively
$(\complex(\{z\}),f(z) \mapsto f(qz) )$ and
$(\complex(\{z^{-1}\}),f(z) \mapsto f(qz) )$). As expected, they
are subgroups of the (global) Galois group. Nevertheless, since
these groups are of second importance in what follows, we omit the
details and
we refer the interesting reader to \cite{sauloyqgaloisfuchs}.\\

\noindent \textbf{Remark.} In \cite{psgaloistheory}, Van der Put
and Singer showed that the Galois groups defined using a
Picard-Vessiot theory can be recovered by means of Tannakian
duality : it is the group of tensor automorphisms of some suitable
complex valued fiber functor over $\mathcal E$. Since two complex
valued fiber functors on a same Tannakian category are necessarily
isomorphic, we conclude that
Sauloy's and Van der Put and Singer's theories coincide.\\

In the rest of this section we exhibit some natural elements of
the Galois group of a given Fuchsian $q$-difference system and we
state the density theorem of Sauloy.

\subsection{The density theorem}

Fix a ``base point" $y_0\in \bon=\complex^* \setminus
\{\text{zeros of } \det(P(z)) \text{ or poles of } P(z)\}$ .
Sauloy exhibits in \cite{sauloyqgaloisfuchs} the following
elements of the (global) Galois group associated to the
$q$-difference system (\ref{general q diff system}) :
\begin{itemize}
\item[1.a)] $\gamma_1(D^{(0)})$ and $\gamma_2(D^{(0)})$ where :
$$\gamma_1:\complex^*=\mathbb{U} \times q^\real \rightarrow \mathbb{U}$$
is the projection over the first factor and : $$\gamma_2 :
\complex^*=\mathbb{U} \times q^\real \rightarrow \complex^*$$ is defined by
$\gamma_2(uq^\omega)=e^{2\pi i \omega}$.
\item[1.b)] $U^{(0)}$.
\item[2.a)] $\breve{P}(y_0)^{-1}\gamma_1(D^{(\infty)})\breve{P}(y_0)$ and $\breve{P}(y_0)^{-1}\gamma_2(D^{(\infty)})\breve{P}(y_0)$.
\item[2.b)] $\breve{P}(y_0)^{-1} U^{(\infty)} \breve{P}(y_0)$.
\item[3)] $\breve{P}(y_0)^{-1}\breve{P}(z)$, $z \in \bon$. \\
\end{itemize}

The following result is due to Sauloy \cite{sauloyqgaloisfuchs}.

\begin{theo}\label{density theo}
The algebraic group generated by the matrices 1.a. to 3. is the
(global) Galois group $G$ of the $q$-difference system
(\ref{general q diff system}). The algebraic group generated by
the matrices 1.a) and 1.b) is the local Galois group at $0$ of the
$q$-difference system (\ref{general q diff system}). The algebraic
group generated by the matrices 2.a) and 2.b) is the local Galois
group at $\infty$ of the $q$-difference system (\ref{general q
diff system}).
\end{theo}

The algebraic group generated by the matrices 3) is called the
\textit{connection component} of the Galois group $G$. The
following result is easy but very useful. Its proof is left to the
reader.

\begin{lem} The connection component of the Galois group $G$ of a regular singular $q$-difference system is a subgroup of the
identity component $G^I$ of $G$.
\end{lem}

\section{Galois groups of the basic hypergeometric equations : non-resonant and non-logarithmic cases}\label{section generique}

We write $a=uq^\alpha$, $b=vq^\beta$ and $c=wq^\gamma$ with
$u,v,w\in\mathbb{U}$ and $\alpha,\beta,\gamma \in \real$ (we
choose a logarithm of $q$).\\

In this section we are aiming to compute the Galois group of the
basic hypergeometric system (\ref{syst hypergeo}) under the
following assumptions :
$$a/b \not \in q^\rinteger \text{ and }c \not \in q^\rinteger.$$

First, we give explicit formulas for the generators of the Galois group of (\ref{syst hypergeo}) involved in Theorem \ref{density theo}.\\

\textit{Local fundamental system of solutions at $0$.} We have :
$$\qhypermatrice{a}{b}{c}{0}=\pmatrice{1}{1}{1}{q/c}\pmatrice{1}{0}{0}{q/c} \pmatrice{1}{1}{1}{q/c}^{-1}.$$
Hence the system (\ref{syst hypergeo}) is non-resonant, and
non-logarithmic at $0$ since $\qhypermatrice{a}{b}{c}{0}$ is
semi-simple. A fundamental system of solutions at $0$ of
(\ref{syst hypergeo}) as described in section \ref{section the
basic objects} is given by
$\yz{a}{b}{c}{z}=\fz{a}{b}{c}{z}e^{(0)}_{\jz{c}}(z)$ with $\jz{c}=
\text{diag}(1,q/c)$ and :
$$\fz{a}{b}{c}{z}=\pmatrice{\qhyper{a}{b}{c}{z}}{\qhyper{aq/c}{bq/c}{q^2/c}{z}}
{\qhyper{a}{b}{c}{qz}}{(q/c) \qhyper{aq/c}{bq/c}{q^2/c}{qz}}.$$\\

\textit{Generators of the local Galois group at $0$.} We have two generators :
$$\pmatrice{1}{0}{0}{e^{2\pi i \gamma}} \text{ and }
\pmatrice{1}{0}{0}{w}.$$

\textit{Local fundamental system of solutions at $\infty$.} We
have :
$$\qhypermatrice{a}{b}{c}{\infty}=\pmatrice{1}{1}{1/a}{1/b} \pmatrice{1/a}{0}{0}{1/b} \pmatrice{1}{1}{1/a}{1/b}^{-1}.$$
Hence the system (\ref{syst hypergeo}) is non-resonant and
non-logarithmic at $\infty$ and a fundamental system of solutions
at $\infty$ of (\ref{syst hypergeo}) as described in section
\ref{section the basic objects} is given by
$\yinf{a}{b}{c}{z}=\finf{a}{b}{c}{z}e^{(\infty)}_{\jinf{a}{b}}(z)$
with $\jinf{a}{b}= \text{diag}(1/a,1/b)$ and :
$$\finf{a}{b}{c}{z}=
\pmatrice{\qhyper{a}{aq/c}{aq/b}{\frac{cq}{ab}z^{-1}}}{\qhyper{b}{bq/c}{bq/a}{\frac{cq}{ab}z^{-1}}}
{\frac{1}{a} \qhyper{a}{aq/c}{aq/b}{\frac{c}{ab}z^{-1}}}{\frac{1}{b} \qhyper{b}{bq/c}{bq/a}{\frac{c}{ab}z^{-1}}}.$$\\

\textit{Generators of the local Galois group at $\infty$.} We have two generators :
$$\breve{P}(y_0)^{-1}\pmatrice{e^{2\pi i \alpha}}{0}{0}{e^{2\pi i \beta}}\breve{P}(y_0) \text{ and } \breve{P}(y_0)^{-1}\pmatrice{u}{0}{0}{v}\breve{P}(y_0).$$\\

\textit{Birkhoff matrix.} The Barnes-Mellin-Watson formula (\textit{cf.} \cite{gasperrahman}) entails that :
$$P(z)= (e^{(\infty)}_{\jinf{a}{b}}(z))^{-1}
\pmatrice{\frac{\pochamer{b,c/a}{q}{\infty}}{\pochamer{c,b/a}{q}{\infty}} \frac{\theta_q (a z)}{\theta_q(z)}}
{\frac{\pochamer{bq/c,q/a}{q}{\infty}}{\pochamer{q^2/c,b/a}{q}{\infty}} \frac{\theta_q (\frac{aq}{c} z)}{\theta_q(z)}}
{\frac{\pochamer{a,c/b}{q}{\infty}}{\pochamer{c,a/b}{q}{\infty}} \frac{\theta_q (b z)}{\theta_q(z)}}
{\frac{\pochamer{aq/c,q/b}{q}{\infty}}{\pochamer{q^2/c,a/b}{q}{\infty}} \frac{\theta_q (\frac{bq}{c} z)}{\theta_q(z)}}
e^{(0)}_{\jz{c}}(z).$$\\

\textit{Twisted Birkhoff matrix.} We have :
$$\breve{P}(z)= \pmatrice{(1/z)^{-\alpha}}{0}{0}{(1/z)^{-\beta}}
\pmatrice{\frac{\pochamer{b,c/a}{q}{\infty}}{\pochamer{c,b/a}{q}{\infty}}
\frac{\theta_q (a z)}{\theta_q(z)}}
{\frac{\pochamer{bq/c,q/a}{q}{\infty}}{\pochamer{q^2/c,b/a}{q}{\infty}}
\frac{\theta_q (\frac{aq}{c} z)}{\theta_q(z)}}
{\frac{\pochamer{a,c/b}{q}{\infty}}{\pochamer{c,a/b}{q}{\infty}}
\frac{\theta_q (b z)}{\theta_q(z)}}
{\frac{\pochamer{aq/c,q/b}{q}{\infty}}{\pochamer{q^2/c,a/b}{q}{\infty}}
\frac{\theta_q (\frac{bq}{c} z)}{\theta_q(z)}}
\pmatrice{1}{0}{0}{z^{1-\gamma}}.$$

$_{}$\vskip 10 pt

We need to consider different cases.

$_{}$\vskip 10 pt

\condun{1} $\underline{a,b,c,a/b,a/c, b/c \not \in q^\rinteger
\text{ and } a/b \text{ or } c \not \in \pm q^{\rinteger/2}}$.

$_{}$\vskip 5 pt

Under this assumption the four numbers $\frac{\pochamer{b,c/a}{q}{\infty}}{\pochamer{c,b/a}{q}{\infty}}$,
$\frac{\pochamer{bq/c,q/a}{q}{\infty}}{\pochamer{q^2/c,b/a}{q}{\infty}}$,
$\frac{\pochamer{a,c/b}{q}{\infty}}{\pochamer{c,a/b}{q}{\infty}}$ and
$\frac{\pochamer{aq/c,q/b}{q}{\infty}}{\pochamer{q^2/c,a/b}{q}{\infty}}$ are non-zero.

\begin{prop}
Suppose that \condun{1} holds. Then the natural action of $G^I$ on $\complex^2$ is irreducible.
\end{prop}

\begin{proof}
Suppose, at the contrary, that the action of $G^I$ is reducible
and let $L \subset \complex^2$ be an invariant line.

Remark that $L$ is distinct from $\complex \vect{1}{0}$ and
$\complex \vect{0}{1}$. Indeed, assume at the contrary that
$L=\complex \vect{1}{0}$ (the case $L=\complex \vect{0}{1}$ is
similar). The line $L=\complex \vect{1}{0}$ being in particular
invariant by the connection component, we see that the line
generated by $\breve{P}(z) \vect{1}{0}$ does not depend on $z \in
\bon$. This yields a contradiction because the ratio of the
components of $\breve{P}(z) \vect{1}{0}=
\vect{\frac{\pochamer{b,c/a}{q}{\infty}}{\pochamer{c,b/a}{q}{\infty}}
\frac{\theta_q (a z)}{\theta_q(z)}(1/z)^{-\alpha}}
{\frac{\pochamer{a,c/b}{q}{\infty}}{\pochamer{c,a/b}{q}{\infty}}
\frac{\theta_q (b z)}{\theta_q(z)}(1/z)^{-\beta}}$ depends on $z$
(remember that $a/b \not \in q^\rinteger$).

On the other hand, since, for all $n\in\integer$, both matrices
$\pmatrice{1}{0} {0}{e^{2\pi i \gamma n}}$ and
$\pmatrice{1}{0}{0}{w^n}$ belong to $G$ and since $G^I$ is a
normal subgroup of $G$, both lines $L_n:=\pmatrice{1}{0}
{0}{e^{2\pi i \gamma n}}L$ and $L'_n:=\pmatrice{1}{0}{0}{w^n} L$
are also invariant by $G^I$.

Note that because \condun{1} holds, at least one of the complex
numbers $w,e^{2\pi i \gamma},u/v,e^{2\pi i (\alpha-\beta)}$ is
distinct from $\pm 1$.

First, suppose that $w \neq \pm 1$. We have seen that $L \neq
\complex \vect{1}{0}, \complex \vect{0}{1}$, hence $L_0,L_1,L_2$
are three distinct lines invariant by the action of $G^I$. This
implies that $G^I$ consists of scalar matrices : this is a
contradiction (because, for instance, $\complex \vect{1}{0}$ is
not invariant for the action of $G^I$). Hence, if $w \neq \pm 1$
we have proved that $G^I$ acts irreducibly.

The case $e^{2\pi i \gamma} \neq \pm 1$ is similar.

Last, the proof is analogous if $u/v \neq \pm 1$ or $e^{2\pi i
(\alpha-\beta)} \neq \pm 1$ (we then use the fact that, for all
$z\in \bon$, $G^I$ is normalized by $\breve{P}(z)^{-1}
\pmatrice{u}{0}{0}{v} \breve{P}(z)$ and $\breve{P}(z)^{-1}
\pmatrice{e^{2\pi i \alpha}}{0}{0}{e^{2\pi i \beta}} \breve{P}(z)$
and that there exists $z \in \bon$ such that $\breve{P}(z)L$ is
distinct from $\complex \vect{1}{0}$ and $\complex \vect{0}{1}$).
\end{proof}

We have the following theorem.

\begin{theo}\label{theo un}
Suppose that \condun{1} holds.
Then we have the following dichotomy :
\begin{itemize}
\item[$\bullet$] if $abq/c \not \in q^\rinteger$ then $G=\Gl_2(\complex)$;
\item[$\bullet$] if $abq/c \in q^\rinteger$ then $G=\overline{\langle \Sl_2(\complex),\sqrt{w}I,e^{\pi i \gamma}I \rangle}$.
\end{itemize}
\end{theo}

\begin{proof}
Since $G^I$ acts irreducibly on $\complex^2$, the general theory
of algebraic groups entails that $G^I$ is generated by its center
$Z(G^I)$ together with its derived subgroup $G^{I,der}$ and that
$Z(G^I)$ acts as scalars. Hence, $G^{I,der} \subset
\Sl_2(\complex)$ also acts irreducibly on $\complex^2$. Therefore
$G^{I,der}=\Sl_2(\complex)$ (a connected algebraic group of
dimension less than or equal to $2$ is solvable hence $\text{dim}
(G^{I,der}) =3$ and $G^{I,der}=\Sl_2(\complex)$). In order to
complete the proof, it is sufficient to determine $\text{det}(G)$.
We have :
\begin{eqnarray*}
\text{det}(\breve{P}(z))
&=&\frac{(1/z)^{-(\alpha+\beta)}z^{1-\gamma}}{\pochamer{q^2/c,a/b,c,b/a}{q}{\infty}}
\left(\underbrace{\theta_q(b)\theta_q(c/a) \frac{\theta_q (a
z)}{\theta_q(z)} \frac{\theta_q (\frac{bq}{c} z)}{\theta_q(z)}
-\theta_q(c/b)\theta_q(a)\frac{\theta_q (\frac{aq}{c}
z)}{\theta_q(z)} \frac{\theta_q (b
z)}{\theta_q(z)}}_{\psi(z)}\right).
\end{eqnarray*}
A straightforward calculation shows that the function :
$$\theta_q(b)\theta_q(c/a) \theta_q (a z) \theta_q (\frac{bq}{c}
z) -\theta_q(c/b)\theta_q(a)\theta_q (\frac{aq}{c} z)\theta_q (b
z)$$ vanishes for $z \in q^\rinteger$ and for $z \in
\frac{c}{abq}q^\rinteger$. On the other hand $\psi$ is a solution
of the first order $q$-difference equation $y(qz)=\frac{c}{abq}
y(z)$. Hence, if we suppose that $abq/c \not \in q^\rinteger$, we
deduce that the ratio $\chi(z)=\frac{\psi(z)}{\frac{\theta_q
(\frac{abq}{c} z)}{\theta_q(z)}}$ defines an holomorphic elliptic
function over $\complex^*$. Therefore $\chi$ is constant and,
evaluating $\chi$ at $z=1/b$,
 we get :
$$\chi=-b \theta_q(a/b)
 \theta_q (c).$$
Finally, we obtain the identity :
\begin{eqnarray} \label{det}
\text{det}(\breve{P}(z))&=&\frac{1-q/c}{1/a-1/b}(1/z)^{-(\alpha+\beta)}z^{1-\gamma}
\frac{\theta_q (\frac{abq}{c} z)}{\theta_q(z)}.
\end{eqnarray}
By analytic continuation (with respect to the parameters) we see
that this formula also holds if $abq/c \in q^\rinteger$.

Consequently, if $abq/c \not \in q^\rinteger$, for any fixed $y_0
\in \bon$, $\text{det}(\breve{P}(y_0)^{-1}\breve{P}(z))$ is a non
constant holomorphic function (with respect to $z$). This implies
that $G=G^I=\Gl_2(\complex)$. On the other hand, if $abq/c \in
q^\rinteger$, then we have that
$\text{det}(\breve{P}(y_0)^{-1}\breve{P}(z))=1$, so that the
connection component of the Galois group is a subgroup of
$\Sl_2(\complex)$ and the Galois group $G$ is the smallest
algebraic group which contains $\Sl_2(\complex)$ and $\{
\sqrt{w}I,e^{\pi i \gamma}I\}$.
\end{proof}

We are going to study the case $a,b,c,a/b,a/c, b/c \not \in q^\rinteger$ and $a/b,c \in \pm q^{\rinteger+1/2}$ in two steps.\\

$_{}$\vskip 10 pt

\condun{2} $\underline{a,b,c,a/b,a/c, b/c \not \in q^\rinteger
\text{ and } q^\rinteger a \cup q^\rinteger b \cup q^\rinteger
aq/c \cup q^\rinteger bq/c=q^\rinteger a \cup -q^\rinteger a \cup
q^{\rinteger+1/2} a \cup -q^{\rinteger+1/2} a}$

$_{}$\vskip 5 pt

 We first establish a preliminary result.

\begin{lem}\label{hg}
Suppose that \condun{2} holds. Then any functional equation of the
form $Az^{n/2}\theta_q(q^Naz)+Bz^{m/2}\theta_q(-q^M
az)+Cz^{l/2}\theta_q(q^Lq^{1/2}az)+Dz^{k/2}
\theta_q(-q^Kq^{1/2}az)=0$ with $A,B,C,D \in \complex$,
$n,m,l,k,N,M,L,K\in \rinteger$ is trivial, that is $A=B=C=D=0$.
\end{lem}

\begin{proof}
Using the non-trivial monodromy of $z^{1/2}$, we reduce the
problem to the case of $n,m,l,k$ being odd numbers. In this case,
using the functional equation $\theta_q(qz)=-z^{-1}\theta_q(z)$,
we can assume without loss of generality that $n=l=m=k=0$. The
expansion of $\theta_q$ as an infinite Laurent series
$\theta_q(z)=\sum_{j\in\rinteger}q^{\frac{j(j-1)}{2}}(-z)^j$
ensures that, for all $j \in \rinteger$, the following equality
holds :
$$A(q^{N})^j+B(-q^{M})^j+C(q^{L+1/2})^j+D(-q^{K+1/2})^j=0.$$ Considering the associated generating series, this implies that :
$$\frac{A}{1-q^Nz}+\frac{B}{1+q^Mz}+\frac{C}{1-q^{L+1/2}z}+\frac{D}{1+q^{K+1/2}z}=0.$$
Hence, considering the poles of this rational fraction, we obtain
$A=B=C=D=0$.
\end{proof}

\begin{prop}
Suppose that \condun{2} holds. Then the natural action of $G^I$ on $\complex^2$ is irreducible.
\end{prop}

\begin{proof}
Suppose, at the contrary, that the action of $G^I$ is reducible
and consider an invariant line $L\subset \complex^2$. In
particular, $L$ is invariant under the action of the connection
component. Consequently, the line $\breve{P}(z)L$ does not depend
on $z\in \bon$. This is impossible using Lemma \ref{hg} (the cases
$L=\complex \vect{1}{0}$ or $\complex \vect{0}{1}$ are excluded by
direct calculation; for the remaining cases consider the ratio of
the coordinates of a generator of $L$ and apply Lemma \ref{hg}).
We get a contraction, hence prove that $G^I$ acts irreducibly.
\end{proof}

\begin{theo}
If \condun{2} holds then we have the following dichotomy :
\begin{itemize}
\item[$\bullet$] if $abq/c \not \in q^\rinteger$ then $G=\Gl_2(\complex)$;
\item[$\bullet$] if $abq/c \in q^\rinteger$ then $G=\overline{\langle \Sl_2(\complex),\sqrt{w}I,e^{\pi i \gamma}I \rangle}$.
\end{itemize}
\end{theo}

\begin{proof}
The proof follows the same lines as that of theorem \ref{theo un}.
\end{proof}

The remaining subcases are $b\in -aq^\rinteger$ and $c\in
-q^\rinteger$; $b\in -aq^{\rinteger+1/2}$ and $c\in
-q^{\rinteger+1/2}$; $b\in
aq^{\rinteger+1/2}$ and $c\in q^{\rinteger+1/2}$.\\

$_{}$\vskip 10 pt

\condun{3}  $\underline{a,b,c,a/b,a/c, b/c \not \in q^\rinteger
\text{ and } b\in -aq^\rinteger \text{ and } c\in -q^\rinteger}$.

$_{}$\vskip 5 pt

 We use the following notations : $b = -aq^\delta
\text{ and } c = -q^\gamma$ with
$\delta=\beta-\alpha,\gamma\in\rinteger$.

The twisted connection matrix takes the following form :
\begin{eqnarray*}
\breve{P}(z) &=&(1/z)^{-\alpha} \pmatrice{
\frac{\pochamer{b,c/a}{q}{\infty}}{\pochamer{c,b/a}{q}{\infty}}
\frac{\theta_q (a z)}{\theta_q(z)}}
 {\frac{\pochamer{bq/c,q/a}{q}{\infty}}{\pochamer{q^2/c,b/a}{q}{\infty}}
 \frac{\theta_q (\frac{aq}{c} z)}{\theta_q(z)}z^{1-\gamma}}
{\frac{\pochamer{a,c/b}{q}{\infty}}{\pochamer{c,a/b}{q}{\infty}}
\frac{\theta_q (b z)}{\theta_q(z)}z^{\delta} }
{\frac{\pochamer{aq/c,q/b}{q}{\infty}}{\pochamer{q^2/c,a/b}{q}{\infty}}
\frac{\theta_q (\frac{bq}{c}
z)}{\theta_q(z)}z^{1+\delta-\gamma}} \\
&=&(1/z)^{-\alpha} \pmatrice{
\frac{\pochamer{b,c/a}{q}{\infty}}{\pochamer{c,b/a}{q}{\infty}}
\frac{\theta_q (a z)}{\theta_q(z)}} {
\frac{\pochamer{bq/c,q/a}{q}{\infty}}{\pochamer{q^2/c,b/a}{q}{\infty}}
 q^{\frac{\gamma(1-\gamma)}{2}}a^{\gamma-1} \frac{\theta_q (-a z)}{\theta_q(z)}}
{\frac{\pochamer{a,c/b}{q}{\infty}}{\pochamer{c,a/b}{q}{\infty}}
q^{\frac{-\delta(\delta-1)}{2}}a^{-\delta} \frac{\theta_q (-a
z)}{\theta_q(z)}}
{\frac{\pochamer{aq/c,q/b}{q}{\infty}}{\pochamer{q^2/c,a/b}{q}{\infty}}
q^{-\frac{(\delta-\gamma +1)(\delta -\gamma)}{2}}(-a)^{\gamma -\delta -1}\frac{\theta_q (a
z)}{\theta_q(z)}}
\end{eqnarray*}

\begin{theo}
Suppose that \condun{3} holds. We have $G=R
\pmatrice{\complex^*}{0}{0}{\complex^*} R^{-1} \cup
\pmatrice{1}{0}{0}{-1} R \pmatrice{\complex^*}{0}{0}{\complex^*}
R^{-1}$ for some $R\in \Gl_2(\complex)$ of the form
$\pmatrice{1}{1}{C}{-C}$, $C\in\complex^*$.
\end{theo}

\begin{proof}
Remark that there exist two nonzero constants $A,B$ such that, for
all $z\in \bon$ :
\begin{eqnarray*}
\breve{P}(-1/a)^{-1}\breve{P}(z)&=&(-a)^\alpha \frac{\theta_q
(-1/a)}{\theta_q(-1)} (1/z)^{-\alpha} \pmatrice{\frac{\theta_q (a
z)}{\theta_q(z)}} {A \frac{\theta_q (-a z)}{\theta_q(z)}} {B
\frac{\theta_q (-a z)}{\theta_q(z)} } {\frac{\theta_q (a
z)}{\theta_q(z)}}\\
&=&(-a)^\alpha \frac{\theta_q (-1/a)}{\theta_q(-1)} (1/z)^{-\alpha} R
\pmatrice{\frac{\theta_q (a z)}{\theta_q(z)} + \sqrt{BA}
\frac{\theta_q (-a z)}{\theta_q(z)}}{0} {0} {\frac{\theta_q (a
z)}{\theta_q(z)} - \sqrt{BA} \frac{\theta_q (-a z)}{\theta_q(z)}}
R^{-1}
\end{eqnarray*}
with $R=\pmatrice{1}{1}{\sqrt{B/A}}{-\sqrt{B/A}}$.

Furthermore, we claim that the functions $X(z):=(1/z)^{-\alpha}
(\frac{\theta_q (a z)}{\theta_q(z)} + \sqrt{BA} \frac{\theta_q (-a
z)}{\theta_q(z)})$ and $Y(z):=(1/z)^{-\alpha} (\frac{\theta_q (a
z)}{\theta_q(z)} - \sqrt{BA} \frac{\theta_q (-a z)}{\theta_q(z)})$
do not satisfy any non-trivial relation of the form $X^rY^s=1$
with $(r,s)\in\rinteger^2 \setminus \{(0,0\}$. Indeed, suppose on
the contrary that such a relation holds. Then
$\frac{((1/z)^{-\alpha} (\theta_q (a z) + \sqrt{BA} \theta_q (-a
z)))^r}{((1/z)^{-\alpha} (\theta_q (a z) - \sqrt{BA} \theta_q (-a
z)))^s}=\theta_q(z)^{s-r}$. Let us first exclude the case $r\neq
s$. If $s>r$ then we conclude that $\theta_q (a z) + \sqrt{BA}
\theta_q (-a z)$ must vanish on $q^\rinteger$. In particular,
$\theta_q (a) + \sqrt{BA} \theta_q (-a)=0$ and $\theta_q (aq) +
\sqrt{BA} \theta_q (-aq) = -(az)^{-1}(\theta_q (a z) - \sqrt{BA}
\theta_q (-a z))=0$, so $\theta_q (a)=0$ that is $a\in
q^\rinteger$. This yields a contradiction. The case $r>s$ is
similar by symmetry. Hence we have $r=s$ so that
$\left(\frac{\theta_q (a z) + \sqrt{BA} \theta_q (-a z)}{\theta_q
(a z) - \sqrt{BA} \theta_q (-a z)}\right)^r=1$. Therefore the
function $\frac{\theta_q (a z) + \sqrt{BA} \theta_q (-a
z)}{\theta_q (a z) - \sqrt{BA} \theta_q (-a z)}$ is constant. This
is clearly impossible and our claim is proved.

This ensures that the connection component of $G^I$, generated by
the matrices $\breve{P}(-1/a)^{-1}\breve{P}(z)$, $z \in \bon$, is
equal to $R\pmatrice{\complex^*}{0}{0}{\complex^*}R^{-1}$.
Consequently, $G$ is generated as an algebraic group by
$R\pmatrice{\complex^*}{0}{0}{\complex^*}R^{-1}$,
$\pmatrice{1}{0}{0}{-1}$,
$\breve{P}(-1/a)^{-1}\pmatrice{u}{0}{0}{-u}\breve{P}(-1/a)=\pmatrice{u}{0}{0}{-u}$
and $\breve{P}(-1/a)^{-1}\pmatrice{e^{2\pi i
\alpha}}{0}{0}{e^{2\pi i \alpha}}\breve{P}(-1/a)=\pmatrice{e^{2\pi
i \alpha}}{0}{0}{e^{2\pi i \alpha}}$. The theorem follows.
\end{proof}

$\bullet$ Both cases ($b\in -aq^{\rinteger+1/2}$ and $c\in -q^{\rinteger+1/2}$) and ($b\in
aq^{\rinteger+1/2}$ and $c\in q^{\rinteger+1/2}$) are similar.\\

$_{}$\vskip 10 pt

\condun{4} $\underline{a \in q^{\integer^*}}$.

$_{}$\vskip 5 pt

In this case, the twisted connection matrix is lower triangular :

\begin{eqnarray*}
\breve{P}(z) &=&\pmatrice{
\frac{\pochamer{b,c/a}{q}{\infty}}{\pochamer{c,b/a}{q}{\infty}}
(-1)^{\alpha} q^{-\frac{\alpha (\alpha -1)}{2}}} {0}
{\frac{\pochamer{a,c/b}{q}{\infty}}{\pochamer{c,a/b}{q}{\infty}}
\frac{\theta_q (b z)}{\theta_q(z)}(1/z)^{-\beta}}
{\frac{\pochamer{aq/c,q/b}{q}{\infty}}{\pochamer{q^2/c,a/b}{q}{\infty}}
\frac{\theta_q (\frac{bq}{c}
z)}{\theta_q(z)}(1/z)^{-\beta}z^{1-\gamma}}
\end{eqnarray*}

\begin{theo}\label{theo theo}
Suppose that \condun{4} holds.
We have the following trichotomy :
\begin{itemize}
\item[$\bullet$] if $b/c \not\in q^\rinteger$ then
$G=\pmatrice{1}{0} {\complex}{\complex^*}$;
\item[$\bullet$] if $c/b \in q^{\integer^*}$ then
$G=\pmatrice{1}{0}{\complex}{\overline{\langle w, e^{2 \pi i
\gamma}\rangle}}$;
\item[$\bullet$] if $bq/c \in q^{\integer^*}$ then
$G=\pmatrice{1}{0} {0}{\overline{\langle w, e^{2 \pi i
\gamma}\rangle}}$.
\end{itemize}
\end{theo}

\begin{proof}
Remark that in each case there exist two constants $A,B$ with $B
\neq 0$ such that, for all $z\in \bon$,
$\breve{P}(1/b)^{-1}\breve{P}(z)= \pmatrice{1}{0}{A\frac{\theta_q
(b z)}{\theta_q(z)}(1/z)^{-\beta}}{B \frac{\theta_q (\frac{bq}{c}
z)}{\theta_q(z)}(1/z)^{-\beta}z^{1-\gamma}}$. Hence, the
connection component is a subgroup of
$\pmatrice{1}{0}{\complex}{\complex^*}$.

Assume $b/c \not\in q^\rinteger$. Then $A\neq 0$ and we claim that
the connection component is equal to
$\pmatrice{1}{0}{\complex}{\complex^*}$. Indeed, for all
$n\in\rinteger$, the following matrix :
$$(\breve{P}(1/b)^{-1}\breve{P}(z))^n=\pmatrice{1}{0}{A\frac{\theta_q
(b z)}{\theta_q(z)}(1/z)^{-\beta}\frac{1-\left(B \frac{\theta_q
(\frac{bq}{c}
z)}{\theta_q(z)}(1/z)^{-\beta}z^{1-\gamma}\right)^n}{1-B
\frac{\theta_q (\frac{bq}{c}
z)}{\theta_q(z)}(1/z)^{-\beta}z^{1-\gamma}}}{\left(B
\frac{\theta_q (\frac{bq}{c}
z)}{\theta_q(z)}(1/z)^{-\beta}z^{1-\gamma}\right)^n}$$ belongs to
the connection component. Consider a polynomial in two variables
$K(X,Y) \in \complex[X,Y]$ such that :
$$K(A\frac{\theta_q (b
z)}{\theta_q(z)}(1/z)^{-\beta}\frac{1-\left(B \frac{\theta_q
(\frac{bq}{c}
z)}{\theta_q(z)}(1/z)^{-\beta}z^{1-\gamma}\right)^n}{1-B
\frac{\theta_q (\frac{bq}{c}
z)}{\theta_q(z)}(1/z)^{-\beta}z^{1-\gamma}},\left(B \frac{\theta_q
(\frac{bq}{c}
z)}{\theta_q(z)}(1/z)^{-\beta}z^{1-\gamma}\right)^n)=0.$$ If $K$
was non zero then we could assume that $K(X,0)\neq 0$. But, for
all $z\in\bon$ in a neighborhood of $\frac{c}{bq}$, we have
$\left|B \frac{\theta_q (\frac{bq}{c}
z)}{\theta_q(z)}(1/z)^{-\beta}z^{1-\gamma}\right|<1$, hence
letting $n$ tend to $+\infty$, we would get $K(A\frac{\theta_q (b
z)}{\theta_q(z)}(1/z)^{-\beta}\frac{1}{1-B \frac{\theta_q
(\frac{bq}{c} z)}{\theta_q(z)}(1/z)^{-\beta}z^{1-\gamma}},0)=0$
which would imply $K(X,0)=0$. This proves that $K=0$. In other
words the only algebraic subvariety of $\complex \times
\complex^*$ containing $(A\frac{\theta_q (b
z)}{\theta_q(z)}(1/z)^{-\beta}\frac{1-\left(B \frac{\theta_q
(\frac{bq}{c}
z)}{\theta_q(z)}(1/z)^{-\beta}z^{1-\gamma}\right)^n}{1-B
\frac{\theta_q (\frac{bq}{c}
z)}{\theta_q(z)}(1/z)^{-\beta}z^{1-\gamma}},\left(B \frac{\theta_q
(\frac{bq}{c}
z)}{\theta_q(z)}(1/z)^{-\beta}z^{1-\gamma}\right)^n)$ for all
$n\in\rinteger$ is $\complex \times \complex^*$ itself. In
particular, the algebraic group generated by the matrix
$(\breve{P}(1/b)^{-1}\breve{P}(z))^n$ for all $n\in\rinteger$ is
$\pmatrice{1}{0}{\complex}{\complex^*}$, hence the connection
component is equal to $\pmatrice{1}{0}{\complex}{\complex^*}$.

It is now straightforward that
$G=\pmatrice{1}{0}{\complex}{\complex^*}$.

Suppose that $c/b \in q^{\integer^*}$. Then the function
$\frac{\theta_q (\frac{bq}{c}
z)}{\theta_q(z)}(1/z)^{-\beta}z^{1-\gamma}$ is constant. Hence the
matrix  $\breve{P}(1/b)^{-1}\breve{P}(z)$ simplifies as follows :
$$\breve{P}(1/b)^{-1}\breve{P}(z)= \pmatrice{1}{0}{A\frac{\theta_q
(b z)}{\theta_q(z)}(1/z)^{-\beta}}{1}$$ with $ A\neq0$. The
connection component is equal to $\pmatrice{1}{0}{\complex}{1}$
and the whole Galois group $G$ is equal to
$\pmatrice{1}{0}{\complex}{\overline{\langle w, e^{2 \pi i
\gamma}\rangle}}$.

Last, suppose that $bq/c \in q^{\integer^*}$. Then $A=0$ and the
function $\frac{\theta_q (\frac{bq}{c}
z)}{\theta_q(z)}(1/z)^{-\beta}z^{1-\gamma}$  is constant, hence
$G=\pmatrice{1}{0}{0}{\overline{\langle w, e^{2 \pi i
\gamma}\rangle}}$.
\end{proof}

$_{}$\vskip 10 pt

\condun{5} $\underline{a \in q^{-\integer}}$.

$_{}$\vskip 5 pt

In this case, the twisted connection matrix is upper triangular :

\begin{eqnarray*}
\breve{P}(x) &=&\pmatrice{
\frac{\pochamer{b,c/a}{q}{\infty}}{\pochamer{c,b/a}{q}{\infty}}
(-1)^{\alpha} q^{-\frac{\alpha (\alpha -1)}{2}}}
{\frac{\pochamer{bq/c,q/a}{q}{\infty}}{\pochamer{q^2/c,b/a}{q}{\infty}}
\frac{\theta_q (\frac{aq}{c}
z)}{\theta_q(z)}(1/z)^{-\alpha}z^{1-\gamma}} {0}
{\frac{\pochamer{aq/c,q/b}{q}{\infty}}{\pochamer{q^2/c,a/b}{q}{\infty}}
\frac{\theta_q (\frac{bq}{c}
z)}{\theta_q(z)}(1/z)^{-\beta}z^{1-\gamma}}
\end{eqnarray*}

\begin{theo}
Suppose that \condun{5} holds.
We have the following trichotomy :
\begin{itemize}
\item[$\bullet$] if $b/c \not\in q^\rinteger$ then
$G=\pmatrice{1}{\complex}{0}{\complex^*}$;
\item[$\bullet$] if $bq/c \in q^{\integer^*}$ then
$G=\pmatrice{1}{\complex}{0}{\overline{\langle w, e^{2 \pi i
\gamma}\rangle}}$;
\item[$\bullet$] if $c/b \in q^{\integer^*}$ then
$G=\pmatrice{1}{0} {0}{\overline{\langle w, e^{2 \pi i
\gamma}\rangle}}$.
\end{itemize}
\end{theo}

\begin{proof}
We argue as for theorem \ref{theo theo}.
\end{proof}

$\bullet$ The cases $\underline{b \in q^\rinteger \text{ or } a/c \in q^\rinteger \text{ or } b/c \in q^\rinteger}$ is similar to the case $a \in q^\rinteger$. We leave the details to the reader.\\


\section{Galois groups of the basic hypergeometric equations : logarithmic cases}

We write $a=uq^\alpha$, $b=vq^\beta$ and $c=wq^\gamma$ with
$u,v,w\in\mathbb{U}$ and $\alpha,\beta,\gamma \in \real$ (we
choose a logarithm of $q$).

\subsection{$c=q$ and $a/b \not \in q^\rinteger$}

The aim of this section is to compute the Galois group of the
basic hypergeometric system (\ref{syst hypergeo}) under the
assumption : $c=q$ and $a/b \not \in q^\rinteger$.

$_{}$

\textit{Local fundamental system of solutions at $0$.} We have :
$$\qhypermatrice{a}{b}{q}{0}=\pmatrice{1}{0}{1}{1} \pmatrice{1}{1}{0}{1} \pmatrice{1}{0}{1}{1}^{-1}.$$
Consequently, we are in the non-resonant logarithmic case at $0$.
We consider this situation as a degenerate case as $c$ tends to $q$, $c\neq q$.

More precisely, we consider the limit as $c$ tends to $q$, with
$c\not \eq q$, of the following matrix-valued function :

\begin{eqnarray*}
&&\fz{a}{b}{c}{z}\pmatrice{1}{1}{1}{q/c}^{-1}\pmatrice{1}{0}{1}{1}\\
&=&\frac{-c}{c-q} \pmatrice{(q/c-1)\qhyper{a}{b}{c}{z}}{\qhyper{aq/c}{bq/c}{q^2/c}{z}-\qhyper{a}{b}{c}{z}}{(q/c-1)\qhyper{a}{b}{c}{qz}}{(q/c) \qhyper{aq/c}{bq/c}{q^2/c}{qz}-\qhyper{a}{b}{c}{qz}} \\
\end{eqnarray*}

Using the notations : $$\qhyperc{a}{b}{z} =
\frac{d}{dc}_{|c=q}\left[\qhyper{a}{b}{c}{z}\right] \text{ and }
\qhyperabc{a}{b}{z} =
\frac{d}{dc}_{|c=q}\left[\qhyper{aq/c}{bq/c}{q^2/c}{z}\right]$$
the above limit is equal to  :
\begin{eqnarray*}
\fz{a}{b}{q}{z}&:=&\pmatrice{\qhyper{a}{b}{q}{z}}{-q(\qhyperabc{a}{b}{z} - \qhyperc{a}{b}{z})}
{\qhyper{a}{b}{q}{qz}}{\qhyper{a}{b}{q}{qz} -q(\qhyperabc{a}{b}{qz} - \qhyperc{a}{b}{qz})}. \\
\end{eqnarray*}
From (\ref{transfo jauge}) we deduce that $\fz{a}{b}{q}{z}$
satisfies
$\fz{a}{b}{q}{qz}\jz{q}=\qhypermatrice{a}{b}{c}{z}\fz{a}{b}{q}{z}$
with $\jz{q}=\pmatrice{1}{1}{0}{1}$. Hence, this matrix being
invertible as a matrix in the field of meromorphic functions, the
matrix-valued function
$\yz{a}{b}{q}{z}=\fz{a}{b}{q}{z}e^{(0)}_{\jz{q}}(z)$ is a
fundamental system of solutions of the
basic hypergeometric equation with $c=q$. Let us recall that $e^{(0)}_{\jz{q}}(z)=\pmatrice{1}{\ell_q(z)}{0}{1}$. \\
\\

\textit{Generators of the local Galois group at $0$.} We have the following generator :
$$\pmatrice{1}{1}{0}{1}.$$\\

\textit{Local fundamental system of solutions at $\infty$.} The
situation is as Section \ref{section generique}. Hence we are in
the non-resonant and non-logarithmic case at $\infty$ and a
fundamental system of solutions at $\infty$ of (\ref{syst
hypergeo}) as described in Section \ref{section the basic objects}
is given by
$\yinf{a}{b}{q}{z}=\finf{a}{b}{q}{z}e^{(\infty)}_{\jinf{a}{b}}(z)$
with $\jinf{a}{b} = \text{diag}(1/a,1/b)$ and :
$$\finf{a}{b}{q}{z}=
\pmatrice{\qhyper{a}{a}{aq/b}{\frac{q^2}{ab}z^{-1}}}{\qhyper{b}{b}{bq/a}{\frac{q^2}{ab}z^{-1}}}
{\frac{1}{a} \qhyper{a}{a}{aq/b}{\frac{q}{ab}z^{-1}}}{\frac{1}{b} \qhyper{b}{b}{bq/a}{\frac{q}{ab}z^{-1}}}.$$\\

\textit{Generators of the local Galois group at $\infty$.} We have two generators :
$$\breve{P}(y_0)^{-1}\pmatrice{e^{2\pi i\alpha}}{0}{0}{e^{2\pi i \beta}}\breve{P}(y_0) \text{ and } \breve{P}(y_0)^{-1}\pmatrice{u}{0}{0}{v}\breve{P}(y_0).$$\\

\textit{Connection matrix}. The connection matrix is the limit as $c$ tends to $q$ of :

\begin{eqnarray*}
(e^{(\infty)}_{\jinf{a}{b}}(z))^{-1}\pmatrice{\frac{\pochamer{b,c/a}{q}{\infty}}{\pochamer{c,b/a}{q}{\infty}}
\frac{\theta_q (a z)}{\theta_q(z)}}
{\frac{\pochamer{bq/c,q/a}{q}{\infty}}{\pochamer{q^2/c,b/a}{q}{\infty}}
\frac{\theta_q (\frac{aq}{c} z)}{\theta_q(z)}}
{\frac{\pochamer{a,c/b}{q}{\infty}}{\pochamer{c,a/b}{q}{\infty}}
\frac{\theta_q (b z)}{\theta_q(z)}}
{\frac{\pochamer{aq/c,q/b}{q}{\infty}}{\pochamer{q^2/c,a/b}{q}{\infty}}
\frac{\theta_q (\frac{bq}{c} z)}{\theta_q(z)}}\pmatrice{1}{1}{1}{q/c}^{-1}\pmatrice{1}{0}{1}{1} e^{(0)}_{\jz{q}}(z)\\
\end{eqnarray*}
which is equal to :

\begin{eqnarray*}
P(z)&:=&(e^{(\infty)}_{\jinf{a}{b}}(z))^{-1} \pmatrice{ \un(a,b;q)
\frac{\theta_q (a z)}{\theta_q(z)}}{
 q(\un_c(a,b;q) -\deux_c(a,b;q))
\frac{\theta_q (a z)}{\theta_q(z)} +az \deux(a,b;q) \frac{\theta'_q
(a z)}{\theta_q(z)}} {\trois(a,b;q) \frac{\theta_q (b
z)}{\theta_q(z)}} { q(\trois_c(a,b;q) -\quatre_c(a,b;q)) \frac{\theta_q (b
z)}{\theta_q(z)} +bz \quatre(a,b;q) \frac{\theta'_q (b
z)}{\theta_q(z)}}
e^{(0)}_{\jz{q}}(z)\\
\end{eqnarray*}
where :
\begin{eqnarray*}
\un(a,b;c)= \frac{\pochamer{b,c/a}{q}{\infty}}{\pochamer{c,b/a}{q}{\infty}}; \ \ \deux(a,b;c)=\frac{\pochamer{bq/c,q/a}{q}{\infty}}{\pochamer{q^2/c,b/a}{q}{\infty}};\\
\trois(a,b;c)=\frac{\pochamer{a,c/b}{q}{\infty}}{\pochamer{c,a/b}{q}{\infty}};
\ \ \quatre(a,b;c)=\frac{\pochamer{aq/c,q/b}{q}{\infty}}{\pochamer{q^2/c,a/b}{q}{\infty}}.\\
\end{eqnarray*}
 and where the subscript $c$ means that we take
the derivative with respect to the third variable $c$.\\

\textit{Twisted connection matrix.}
\begin{small} \begin{eqnarray*}
\breve{P}(z)&=&\pmatrice{(1/z)^{-\alpha}}{0}{0}{(1/z)^{-\beta}}
\pmatrice{ \un(a,b;q) \frac{\theta_q (a z)}{\theta_q(z)}}{
 q(\un_c(a,b;q) -\deux_c(a,b;q))
\frac{\theta_q (a z)}{\theta_q(z)} +az \deux(a,b;q)
\frac{\theta'_q (a z)}{\theta_q(z)}} {\trois(a,b;q) \frac{\theta_q
(b z)}{\theta_q(z)}} { q(\trois_c(a,b;q) -\quatre_c(a,b;q))
\frac{\theta_q (b z)}{\theta_q(z)} +bz \quatre(a,b;q)
\frac{\theta'_q (b z)}{\theta_q(z)}} \pmatrice{1}{\ell_q(z)}{0}{1}
\end{eqnarray*}\end{small}

$_{}$\vskip 10 pt

We need to consider different cases.

$_{}$\vskip 10 pt

\conddeux{1} $\underline{a \not\in q^\rinteger \text{ and } b\not \in
q^\rinteger}$.

$_{}$\vskip 5 pt

Subject to this condition, the complex numbers $\un(a,b;q)$,
$\deux(a,b;q)$, $\trois(a,b;q)$ and $\quatre(a,b;q)$ are non-zero.

\begin{prop} \label{blabla}
If \conddeux{1} holds then the natural action of $G^I$ on $\complex^2$ is irreducible.
\end{prop}

\begin{proof} Assume, at the contrary, that the action of $G^I$ is
reducible and let $L$ be an invariant line.

Let us fist remark that $L \neq \complex \vect{1}{0}$ (in
particular, $G^I$ does not consist of scalar matrices). Indeed, if not,
$\complex \vect{1}{0}$ would be stabilized by the connection
component and the line spanned by $\breve{P}(z) \vect{1}{0}$
would be independent of $z\in\bon$ : this is clearly false.

The group $G^I$ being normalized by $\pmatrice{1}{1}{0}{1}$
(since $G^I$ is a normal subgroup of $G$), the
lines $\pmatrice{1}{1}{0}{1}^n L$ are also invariant by the action
of $G^I$. These lines being distinct (since $L\neq \complex
\vect{1}{0}$) we conclude that $G^I$ consists of scalar matrices and
we get a contradiction. This proves that $G^I$ acts
irreducibly.
\end{proof}

As a consequence we have the following theorem.

\begin{theo}
Suppose that \conddeux{1} holds. Then we have the following dichotomy :
\begin{itemize}
\item[$\bullet$] if $ab \not \in q^\rinteger$ then $G=\Gl_2(\complex)$;
\item[$\bullet$] if $ab \in q^\rinteger$ then $G= \Sl_2(\complex)$.
\end{itemize}
\end{theo}

\begin{proof}
Using the irreducibility of the natural action of $G^I$ and
arguing as for the proof of theorem \ref{theo un}, we obtain the
equality $G^{I,der}=\Sl_2(\complex)$. From the formula (\ref{det})
we deduce that the determinant of the twisted connection matrices
when $c=q$ is equal to the limit as $c$ tends to $q$ of
$\frac{-1}{1/a-1/b}(1/z)^{-(\alpha+\beta)} z^{1-\gamma}
\frac{\theta_q (\frac{abq}{c} z)}{\theta_q(z)}$ \textit{i.e.}
$\frac{-1}{1/a-1/b}(1/z)^{-(\alpha+\beta)} z^{1-\gamma}
\frac{\theta_q (ab z)}{\theta_q(z)}.$

If $ab \not \in q^\rinteger$ then this determinant is a
non-constant holomorphic function and consequently
$G=\Gl_2(\complex)$.

If $ab \in q^\rinteger$ then this determinant does not depend on
$z$. This implies that the connection component of the Galois
group is a sub-group of $\Sl_2(\complex)$.  Furthermore, $ab \in
q^\rinteger$ entails that $uv =1$ and $\alpha+\beta \in
\rinteger$, that is, $e^{2\pi i (\alpha+\beta)}=1$. Consequently,
the local Galois groups are subgroups of $\Sl_2(\complex)$ and the
global Galois group $G$ is therefore a subgroup of
$\Sl_2(\complex)$.
\end{proof}



$_{}$\vskip 10 pt

\conddeux{2} $\underline{b\in q^{\integer^*}}$.

$_{}$\vskip 5 pt

Then the twisted
connection matrix simplifies as follows :
\begin{small}\begin{eqnarray*}
\breve{P}(z)&=& \pmatrice{ \un(a,b;q) \frac{\theta_q (a
z)}{\theta_q(z)}(1/z)^{-\alpha}} {q(\un_c(a,b;q) -\deux_c(a,b;q))
\frac{\theta_q (a z)}{\theta_q(z)}(1/z)^{-\alpha} +az \deux(a,b;q)
\frac{\theta'_q (a z)}{\theta_q(z)}(1/z)^{-\alpha}} {0} {
q(\trois_c(a,b;q) -\quatre_c(a,b;q)) (-1)^\beta
q^{-\frac{\beta(\beta-1)}{2}}}
\pmatrice{1}{\ell_q(z)}{0}{1}\\
\end{eqnarray*}\end{small}

\begin{theo}\label{theou}
Suppose that \conddeux{2} holds.
Then we have $G=\pmatrice{\complex^*}{\complex}{0}{1}$.
\end{theo}
\begin{proof}
Fix a point $y_0 \in \bon$ such that $\breve{P}(y_0)$ is of the
form $\pmatrice{A}{B}{0}{C}$ with $A,C \neq 0$. There exists a
constant $D\in \complex^*$ such that :
$$\breve{P}(y_0)^{-1}\breve{P}(z)= \pmatrice{D\frac{\theta_q (a z)}{\theta_q(z)}(1/z)^{-\alpha}}{*}{0}{1}.$$
Since $G^I$ is normalized by $\pmatrice{1}{1}{0}{1}$
(remember that $G^I$ is a normal subgroup of $G$), it contains, for all $n\in \rinteger$, the matrix :
$$\pmatrice{D\frac{\theta_q (a z)}{\theta_q(z)}(1/z)^{-\alpha}}{*+n(D\frac{\theta_q (a z)}{\theta_q(z)}(1/z)^{-\alpha}-1)}{0}{1}.$$
Because $a\not \in q^\rinteger$, the function $D\frac{\theta_q (a
z)}{\theta_q(z)}(1/z)^{-\alpha}-1$ is not identically equal to
zero over $\complex^*$ and therefore $G^I$ contains, for all $z
\in \bon$ :
$$\pmatrice{D\frac{\theta_q (a z)}{\theta_q(z)}(1/z)^{-\alpha}}{\complex}{0}{1}.$$
In particular, $\pmatrice{D\frac{\theta_q (a
z)}{\theta_q(z)}(1/z)^{-\alpha}}{0}{0}{1}$ belongs to $G^I$, so that
$\pmatrice{\complex^*}{0}{0}{1}$ is a subgroup of $G^I$ and
$\pmatrice{\complex^*}{\complex}{0}{1} \subset G$. The converse
inclusion is clear.
\end{proof}

$_{}$\vskip 10 pt

\conddeux{3} $\underline{b\in q^{-\integer}}$.

$_{}$\vskip 5 pt

Using the identity :
$$bz\frac{\theta_q'(bz)}{\theta_q(z)}=(-\beta-\ell_q(z))(-1)^\beta
q^{-\frac{\beta(\beta-1)}{2}}$$ we see that, in this case, the
twisted connection matrix takes the form :

\begin{eqnarray*}
\breve{P}(z)
&=&\pmatrice{0} {q(\un_c(a,b;q) -\deux_c(a,b;q)) \frac{\theta_q (a
z)}{\theta_q(z)}(1/z)^{-\alpha}} {w(a,b,q) (-1)^\beta
q^{-\frac{\beta(\beta-1)}{2}}} { q(\trois_c(a,b;q)
-\quatre_c(a,b;q)-\beta/q) (-1)^\beta q^{-\frac{\beta(\beta-1)}{2}}}\\
\end{eqnarray*}

\begin{theo}
Suppose that \conddeux{3} holds. Then we have $G=\pmatrice{1}{\complex}{0}{\complex^*}.$
\end{theo}

\begin{proof}
Fix a base point $y_0 \in \bon$. There exist three constants
$C,C',C''\in \complex$ with $C\neq 0$, such that the following
identity holds for all $z\in \bon$ :
$$\breve{P}(y_0)^{-1}\breve{P}(z)=\pmatrice{1}{C'\frac{\theta_q (a z)}{\theta_q(z)}(1/z)^{-\alpha}+C''}
{0}{C\frac{\theta_q (a z)}{\theta_q(z)}(1/z)^{-\alpha}}.$$ The proof is similar to the proof of Theorem \ref{theou}.
\end{proof}

The remaining case $a\in q^\rinteger$ is similar to
\conddeux{3}.

The case $a=b$ and $c\not \in q^\rinteger$ is similar to the case treated in this section.

\subsection{$a=b$ and $c=q$}

The aim of this section is to compute the Galois group of the
basic hypergeometric system (\ref{syst hypergeo}) under the
assumption : $a=b$ and $c=q$.

$_{}$

\textit{Local fundamental system of solutions at $0$.} The
situation is the same
as in the case $c=q$ and $a/b \not\in q^\rinteger$.\\

\textit{Generator of the local Galois group at $0$}. We have the following generator :
$$\pmatrice{1}{1}{0}{1}.$$\\

\textit{Local fundamental system of solutions at $\infty$.} We
have :
$$\qhypermatrice{a}{a}{q}{z}=\pmatrice{1}{0}{1/a}{1} \pmatrice{1/a}{1}{0}{1/a} \pmatrice{1}{0}{1/a}{1}^{-1}.$$
Consequently, we are in the non-resonant logarithmic case at
$\infty$. We consider the case $a=b$ and $c=q$ as a degenerate
case of the situation $c=q$ as $a$ tends to $b$, $a/b \neq 1$.

We consider the following matrix valued function :
$$\finf{a}{b}{q}{z}\pmatrice{1}{1}{1/a}{1/b}^{-1}\pmatrice{1}{0}{1/a}{1}.$$
A straightforward calculation, which we omit here because it is
long although easy, shows that this matrix-valued function does admit a limit as $a$ tends to $b$ that we denote $\finf{a}{a}{q}{z}$. A fundamental
system of solutions at $\infty$ of (\ref{syst hypergeo}) as
described in section \ref{section the basic objects} is given by
$\yinf{a}{a}{q}{z}=\finf{a}{a}{q}{z}e^{(\infty)}_{\jinf{a}{a}}(z)$
with $\jinf{a}{a} = \pmatrice{1/a}{1}{0}{1/a}$.\\
\\

\textit{Generator of the local Galois group at $\infty$}. We have the following generators :
$$\pmatrice{u}{0}{0}{u}, \ \ \pmatrice{e^{2\pi i \alpha}}{0}{0}{e^{2\pi i \alpha}} \text{ and } \breve{P}(y_0)^{-1}\pmatrice{1}{a}{0}{1}\breve{P}(y_0).$$\\

\textit{Birkhoff matrix}. The Birkhoff matrix is equal to
$(e^{(\infty)}_{\jinf{a}{a}}(z))^{-1}Qe^{(0)}_{\jz{q}}(z)$ where
$Q$ is the limit as $a$ tends to $b$ of :

\begin{eqnarray*}
&&\pmatrice{1}{0}{1/a}{1}^{-1}\pmatrice{1}{1}{1/a}{1/b}\pmatrice{
\un(a,b;q) \frac{\theta_q (a z)}{\theta_q(z)}}{
 q(\un_c(a,b;q) -\deux_c(a,b;q))
\frac{\theta_q (a z)}{\theta_q(z)} +az \deux(a,b;q)
\frac{\theta'_q (a z)}{\theta_q(z)}} {\trois(a,b;q) \frac{\theta_q
(b z)}{\theta_q(z)}} { q(\trois_c(a,b;q) -\quatre_c(a,b;q))
\frac{\theta_q (b z)}{\theta_q(z)} +bz \quatre(a,b;q)
\frac{\theta'_q (b z)}{\theta_q(z)}}.
\end{eqnarray*}
It has the following form :


\begin{eqnarray*}
&& \pmatrice{ C \frac{\theta_q(az)}{\theta_q(z)} +az
\frac{\theta_q(a)}{\pochamer{q}{q}{\infty}^2} \frac{\theta_q'(a
z)}{\theta_q(z)}} {*} {-(1/a)
\frac{\theta_q(a)}{\pochamer{q}{q}{\infty}^2}\frac{\theta_q(az)}{\theta_q(z)}}
{C'\frac{\theta_q(az)}{\theta_q(z)}-z\frac{\theta_q(a)}{\pochamer{q}{q}{\infty}^2}
\frac{\theta_q'(a z)}{\theta_q(z)}}
\end{eqnarray*}
where $*$ denotes some meromorphic function.\\

\textit{Twisted Birkhoff matrix}.

$$
 (1/z)^{-\alpha} \pmatrice{1}{-a\ell_q(z)}{0}{1}
\pmatrice{C \frac{\theta_q(az)}{\theta_q(z)}+az\frac{\theta_q(a)}{\pochamer{q}{q}{\infty}^2}
\frac{\theta_q'(a z)}{\theta_q(z)}} {*} {-(1/a)
\frac{\theta_q(a)}{\pochamer{q}{q}{\infty}^2}\frac{\theta_q(az)}{\theta_q(z)}}
{C'\frac{\theta_q(az)}{\theta_q(z)}-z\frac{\theta_q(a)}{\pochamer{q}{q}{\infty}^2}
\frac{\theta_q'(a z)}{\theta_q(z)}}
 \pmatrice{1}{\ell_q(z)}{0}{1}.
$$

$_{}$\vskip 10 pt

We need to consider different cases.

$_{}$\vskip 10 pt

\condtrois{1} $\underline{a \not \in q^\rinteger}$.

$_{}$\vskip 5 pt

\begin{prop}
Suppose that \condtrois{1} holds. Then the natural action of $G^I$ on $\complex^2$ is irreducible.
\end{prop}
\begin{proof}
Remark that $\complex \vect{1}{0}$ is not an invariant line.
Indeed, if not, this line would be invariant by the action of the
connection component, hence the line spanned by $\breve{P}(z)
\vect{1}{0}$ would be independent of $z\in\bon$. Considering the
ratio of the coordinates of this line, this would imply the
existence of some constant $A \in \complex$ such that the
following functional equation holds on $\complex^*$ :
$$C \frac{\theta_q(az)}{\theta_q(z)} +az
\frac{\theta_q(a)}{\pochamer{q}{q}{\infty}^2} \frac{\theta_q'(a
z)}{\theta_q(z)}+
\frac{\theta_q(a)}{\pochamer{q}{q}{\infty}^2}\frac{\theta_q(az)}{\theta_q(z)}
\ell_q(z) = A \frac{\theta_q(az)}{\theta_q(z)}.$$ The fact that
$\theta_q(az)$ vanishes exactly to the order one at $z=1/a$,
yields a contradiction.

The end of the proof is similar to the proof of Proposition
\ref{blabla}.
\end{proof}

\begin{theo}
If \condtrois{1} holds then we have the following dichotomy :
\begin{itemize}
  \item[$\bullet$] if $a^2 \not \in q^\rinteger$ then
  $G=\Gl_2(\complex)$;
  \item[$\bullet$] if $a^2 \in q^\rinteger$ then
  $G=\Sl_2(\complex)$.
\end{itemize}
\end{theo}

\begin{proof}
The proof follows the same line as that of theorem \ref{theo un}.
\end{proof}

$_{}$\vskip 10 pt

\condtrois{2} $\underline{a \in q^\rinteger}$.

$_{}$\vskip 5 pt

Under this condition, the connection matrix simplifies as follows,
for some constants $C,C'\in \complex$ :
\begin{eqnarray*}
&&
\pmatrice{ C  } {*}
{0}{C'}.\\
\end{eqnarray*}

\begin{theo}
Suppose that \condtrois{2} holds. Then we have : $G = \pmatrice{1}{\complex}{0}{1}$.
\end{theo}
\begin{proof}
The local Galois group at $0$ is generated by
$\pmatrice{1}{1}{0}{1}$, hence $G$ contains
$\pmatrice{1}{\complex}{0}{1}$.

Since the twisted connection matrix is upper triangular with
constant diagonal entries, the connection component is a subgroup of
$\pmatrice{1}{\complex}{0}{1}$. The generators of the local Galois
group at $0$ and at $\infty$ also lie in
$\pmatrice{1}{\complex}{0}{1}$. Therefore, $G$ is a subgroup of
$\pmatrice{1}{\complex}{0}{1}$.
\end{proof}

\bibliography{biblio}
\bibliographystyle{plain}

\noindent \textsc{Julien Roques\\
D\'epartement de Math\'ematiques et Applications (DMA)\\
\'Ecole Normale
Sup\'erieure\\
45, rue d'Ulm\\
F 75230 Paris cedex 05}\\
E-mail : \textsf{julien.roques@ens.fr}

\end{document}